\documentclass[12pt]{article}
\usepackage{amsmath}
\usepackage{amssymb}
\newcommand{\Z}{\mathbb Z}

\newcommand{\Q}{\mathbb Q}
\newcommand{\N}{\mathbb N}
\newcommand{\C}{\mathbb C}
\newcommand{\R}{\mathbb R}
\newcommand{\p}{\partial}
\newcommand{\z}{\zeta}
\newcommand{\Deg}{{\rm Deg}}
\begin{document}

\title{The Amazing Image Conjecture}
\author{Arno van den Essen}
\maketitle

\begin{abstract}

\noindent In this paper we discuss a general framework in which we present
a new conjecture, due to Wenhua Zhao, the Image Conjecture. This conjecture implies the
Generalized Vanishing Conjecture and hence the Jacobian Conjecture.
 Crucial ingredient is the notion of a Mathieu space: let $k$ be a field
and $R$ a $k$-algebra. A $k$-linear subspace $M$ of $R$ is called a Mathieu
subspace of $R$, if the following holds: let $f\in R$ be such that $f^m\in M$, for
all $m\geq 1$, then for every $g\in R$ also $gf^m\in M$, for almost all $m$,
i.e.\@ only finitely many exceptions.

 Let $A$ be the polynomial ring in $\z=\z_1,\ldots,\z_n$ and $z_1,\ldots,z_n$
over $\C$. The Image Conjecture (IC) asserts that $\sum_i(\partial_{z_i}-\z_i)A$
is a Mathieu subspace of $A$. We prove this conjecture for $n=1$. Also we
relate (IC) to the following Integral Conjecture: if $B$ is an open subset
of $\R^n$ and $\sigma$ a positive measure, such that the integral over $B$
of each polynomial in $z$ over $\C$ is finite, then the set of polynomials,
whose integral over $B$ is zero, is a Mathieu subspace of $\C[z]$.
It turns out that Laguerre polynomials play a special role in the study
of the Jacobian Conjecture. 
\end{abstract}

\section*{Introduction}

Some twenty five years ago I learnt about
the existence of the Jacobian Conjecture, during one of
my visits to my friend Pascal Adjamagbo in Paris. 
The problem always stayed somewhere in my mind and in the
meantime I worked on different related topics and found
counter examples to various conjectures, which would have implied 
the truth of the still mysterious Jacobian Conjecture. All these 
experiences fed my believe that the Jacobian Conjecture, if true 
at all, would be difficult to generalize, since it felt like a 
kind of optimal statement. Therefore I often stated in public the
following dictum

\smallskip
{\leftskip=0.5cm{\rightskip=0.5cm
\noindent {\em ``If you have a conjecture which implies the 
Jacobian Conjecture, but is not equivalent to it, then you 
can be sure that your conjecture is false."}\par}\par}
\smallskip

\noindent It is therefore no surprise that, when in July 2009 Wenhua Zhao came up
with a new conjecture implying the Jacobian Conjecture, I set out to
find a counter example. This conjecture, which was given the name
{\em Image Conjecture} by its inventor, is so general that I was convinced
that it would be easy to find a counterexample. Surprisingly, I did not. Instead I found various instances in favour of it.

The aim of this paper is to bring this fascinating new conjecture to
the attention of a larger audience. Hopefully it will inspire the reader
to join me in my search for either a proof, or counterexample.
The style in which it is written will be easy going. Sometimes
I will skip proofs and refer to the papers of Zhao and the upcoming
joint work with Wright and Zhao ([EWrZ]).

\section{The Image Conjecture: a first encounter}

To please those readers who cannot wait to see what the Image
Conjecture is all about, I will start this section, by describing
its most important special case.

Let $k$ be any field, $A$ a commutative $k$-algebra and $A[z]$ the polynomial
ring in $n$ variables $z=(z_1,\ldots,z_n)$ over $A$. Elements of the ring
$A[z]$ will simply be called polynomials, without refering to $A$ or $z$.
Let $a_1,\ldots,a_n$ be elements of $A$ and denote by $D$ the set of
commuting differential operators
$$\partial_{z_1}-a_1,\ldots,\partial_{z_n}-a_n.$$
\noindent Finally put
$$ImD=\sum_{i=1}^n (\partial_{z_i}-a_i)A[z].$$

\medskip

\noindent {\bf Image Conjecture (IC(n,A)).} {\em Assume that $(a_1,\ldots,a_n)$
is a regular sequence in $A$. If all positive powers
of a polynomial $f$ belong to $ImD$, then for every polynomial $g$,
almost all polynomials $gf^m$ also belong to $ImD$}.

\medskip

The sentence {\em almost all} means for all, with only
a finite number of exceptions. Furthermore the sequence $(a_1,\ldots,a_n)$ is
called a {\em regular sequence in $A$} if $a_1$ is no zero-divisor in $A$,
for each $i\geq 1$ the element $a_{i+1}$ is no zero-divisor
in $A/(a_1,\ldots,a_i)$ and the ideal generated by all $a_i$ is not equal
to $A$.

To get a feeling for the difficulty of the problem, the reader is 
invited to find a proof for the one dimensional case. In fact,
in this dimension the conjecture has only been proved in case the
ideal $Aa_1$ is a radical ideal. If additionaly $A$ is a UFD, also the non-radical
case has been proved.
As we will see below, the Jacobian Conjecture follows from the very
special case where $A$ is the polynomial ring
$\C[\zeta_1,\ldots,\zeta_n]$ and $a_i=\zeta_i$ for each $i$.

The property concerning the powers of polynomials, which is used in the
formulation of the Image Conjecture above, was formalized by Zhao in
his paper [Z4] as follows.

 Let $k$ be a field, $R$ a $k$-algebra (not necessarily
commutative) and $M$ a $k$-linear subspace of $R$. An element
$f$ of $R$ is said to have the {\em left Mathieu property with respect to $M$},
if the following holds: if all positive powers of $f$ belong to $M$,
then for every $g$ of $R$, almost all elements $gf^m$ belong to $M$.
Furthermore $M$ is called a {\em left Mathieu subspace of $R$}, if all
elements of $R$ have the left Mathieu property, with respect to $M$. Similarly
one defines the notion of a right Mathieu subspace and finally
$M$ is called a Mathieu subspace, if it is both a left and right Mathieu
subspace. In most examples discussed in this paper, the ring $R$ will
be commutative, so we just speak about Mathieu subspaces.

\medskip

\noindent Using this terminology, the Image Conjecture formulated above
simply states that $ImD$ is a Mathieu subspace of $A[z]$.

The notion of Mathieu subspace was first introduced by Zhao in [Z4] and
was inspired by the following conjecture proposed by Olivier Mathieu
in [M], 1995.

\medskip

\noindent {\bf Mathieu Conjecture.} {\em Let $G$ be a compact connected real Lie
group with Haar measure $\sigma$. Let $f$ be a complex valued $G$-finite
function on $G$, such that $\int_G f^m d\sigma=0$ for all positive $m$.
Then for every $G$-finite function $g$ on $G$, also $\int_G gf^m d\sigma=0$,
for almost all $m$.}

\medskip

\noindent Here a function $f$ is called $G$-{\em finite}, if the $\C$-vector space
generated by the elements of the orbit $G\cdot f$ is finite dimensional.

With the terminology introduced above, the Mathieu Conjecture
can be reformulated as follows: let $R$ be the $\C$-algebra of complex
valued $G$-finite functions on $G$. Then the $\C$-subspace of $f$'s,
which satisfy $\int _G f d\sigma=0$, is a Mathieu subspace of $R$.

The importance of Mathieu's conjecture comes from the fact that it
implies the Jacobian conjecture, as was shown in [M]. Since its
formulation, only one non-trivial case of this conjecture was solved, namely
the case that $G$ is commutative. This result, which is due to Duistermaat
and van der Kallen, can be formulated as follows (see [DvK]).

\medskip

\noindent {\bf Duistermaat-van der Kallen theorem.} {\em Let $k$ be a field of
characteristic zero and 
$R=k[z_1,\ldots,z_n,z_1^{-1},\ldots,z_n^{-1}]$, the ring of Laurent polynomials over
$k$. Then the set of Laurent polynomials, whose constant term is zero, is
a Mathieu subspace of $R$.}

\medskip

\noindent Already the proof of the one dimensional case is not at all obvious
and again the reader is invited to find an elementary proof.
The hypothesis that $k$ has characteristic zero, cannot be dropped,
as can be seen from the following example, which is due to Roel Willems.

\medskip

\noindent {\bf Counterexample 1.1.} {\em Let $n=1$ and write $t$ instead of $z_1$.
Let $f=t^{-1}+t^{p-1}\in k[t,t^{-1}]$, where $k$ is a field of characteristic
$p>0$. Then the constant term of all positive powers of
$f$ is zero, however for all $m=p^k-1$, the constant coefficient of
$t^{-1}f^m$ is non-zero.}

\medskip

\indent The notion of Mathieu subspaces of a ring $R$ can be viewed as a 
{\em generalization of ideals rings}, since obviously an ideal
of $R$ is a Mathieu subspace of $R$. However, Mathieu subspaces
are far more complicated to understand and to recognize. For example
it is easy to describe all ideals of the univariate polynomial ring
$k[t]$, but even for concrete cases we have no way, other than
ad hoc methods, to decide if a given $k$-linear subspace of $k[t]$ is
a Mathieu subspace or not. The reader who wants to test this statement
is refered to section six, where we discuss various Mathieu subspaces
of $k[t]$ and some candidate ones.

There is however one easy property that Mathieu subspaces share with
ideals and which often can be used to show that a given subspace is
not a Mathieu space.

\medskip

\noindent {\bf The $1$-property.} {\em Let $M$ be a Mathieu subspace of a
$k$-algebra $R$. If $M$ contains $1$, then $M=R$.}

\medskip

\noindent Indeed, if $1$ belongs to $M$, then all positive powers of $1$
belong to $M$. Hence, by the Mathieu property, it follows that for each
$g$ in $R$ almost all elements $g\cdot 1^m$ belong to $M$,\@ i.e. each such
$g$ belongs to $M$. So $M=R$.

\medskip

\noindent {\bf Example 1.2.} {\em Let $R=k[z]$ be the univariate polynomial ring and
$k$ a field of positive characteristic $p$. Let $D=\partial_z$.
Then $ImD$ is not a Mathieu subspace of $R$.}

\medskip

\noindent Namely obviously $ImD$ contains 1, but it does not contain
$z^{p-1}$ since $p=0$ in $R$. So by the $1$-property $ImD$ cannot be a Mathieu
subspace of $R$.

This example shows that the hypothesis concerning the regularity of
the sequence of $a_i$'s in the statement of the Image Conjecture cannot
be dropped.

\section{Motivation for the Image Conjecture}

To understand where the Image Conjecture comes from, we recall some
recent results concerning the Jacobian Conjecture.

As is well-known the Jacobian Conjecture was formulated by O. Keller
in 1939 in his paper [K]. It asserts that if the jacobian determinant
of a polynomial map from complex affine $n$-space to itself is a
non-zero constant, then the map is invertible, in the sense that its
inverse is again a polynomial map. The conjecture is open in all
dimensions $n$ greater than one. In 1982, Bass, Connell and Wright,
and independently Yagzhev showed that in order to prove or disprove
the conjecture it suffices to study
so-called {\em cubic homogeneous} polynomial maps i.e. maps of the form

$$z+H=(z_1+H_1,\ldots,z_n+H_n)$$

\noindent where the $H_i$ are either zero or homogeneous of degree three.
This result is known as the {\em cubic homogeneous reduction}. It
was also shown that the condition for the jacobian determinant of such a
map to be a non-zero constant is equivalent to the nilpotency of the
jacobian matrix of $H$ (see [BCW] or [E]).

In 2003 Michiel de Bondt and the author improved upon the above reduction
result, by showing that one may {\em additionaly assume} that the Jacobian
matrix of $H$ is symmetric (see [BE]), which by Poincar\'e's lemma
implies that $H$ equals the gradient of some (quartic) polynomial
$P$ in $n$ variables over $\C$ i.e. $P\in\C[z]$. Using this fact
Wenhua Zhao obtained the following new and surprising
{\em equivalent} description of the Jacobian Conjecture ([Z1], 2004)

\medskip

\noindent {\bf Vanishing Conjecture.} {\em Let $\Delta=\sum_i {\partial^2_{z_i}}$
be the Laplace operator and let $P\in\C[z]$ be homogeneous. 
If $\Delta^m(P^m)=0$ for all positive $m$, then $\Delta^m(P^{m+1})=0$ 
for almost all $m$.}

\medskip

\noindent In fact the condition $\Delta^m(P^m)=0$ for all
positive $m$ is equivalent to the nilpotency of
the jacobian matrix $J(\nabla(P))$ and the vanishing of
all sufficiently large powers $\Delta^m(P^{m+1})$ is equivalent
to the invertibility of the map  $z+\nabla(P)$. 

Zhao observed the resemblance with Mathieu's conjecture
and could make this resemblance even better by showing that the
Vanishing Conjecture is equivalent to the following version (see [EZ])

\medskip

\noindent {\bf Vanishing Conjecture.} {\em If $P\in\C[z]$ is homogeneous and 
such that $\Delta^m(P^m)=0$ for all positive $m$, then for each $Q$ in 
$\C[z]$ $\Delta^m(QP^m)=0$ for almost all $m$.}

\medskip

\noindent Now it is not difficult to show that if the Vanishing Conjecture
holds for the Laplace operator, it also holds for all quadratic homogeneous
operators with constant coefficients (use Lefschetz's principle and the fact 
that over the complex numbers all quadratic forms are essentially sums of 
squares. For more details we refer to [Z2]).

After these observations Zhao made the following more general conjecture,
dropping the homogeneity condition on the polynomial $P$ and replacing
the Laplace operator by any differential operator with {\em constant
coefficients}

\medskip

\noindent {\bf Generalized Vanishing Conjecture(GVC(n)).} {\em Let $\Lambda$
be any differential operator with constant coefficients, i.e.\@
$\Lambda\in\C[\partial_1,\ldots,\partial_n]$. If $P\in\C[z]$ is such that
$\Lambda^m(P^m)=0$ for all positive $m$, then for each $Q$ in $\C[z]$
also $\Lambda^m(QP^m)=0$ for almost all $m$.}

\medskip

\noindent When I saw this conjecture I was convinced that it should
be easy to find a counterexample.
So I first investigated the one dimensional case. Then the conjecture
is easily seen to be true, namely let $\Lambda$ be a polynomial
in $\C[\partial]$ of order $r\geq 0$, i.e.\@ $\partial^r$
is the lowest degree monomial in $\partial$ appearing in $\Lambda$. Now
let $P$ in $\C[z]$ be a polynomial of degree $d$. Observe that the
polynomial $\Lambda(P)$  has degree $d-r$ if $r\leq d$. In particular,
if $\Lambda(P)=0$ it follows that $r\geq d+1$. Consequently the order
of $\Lambda^m$, which equals $rm$, is greater or equal to $dm+m$, which is 
greater than the degree of $QP^m$ if $m$ is greater than the degree of $Q$. 
This implies that $\Lambda^m(QP^m)=0$ for such $m$.

Next I investigated the special two variable case $\partial_1^p+\partial_2^q$
where $p$ and $q$ are natural numbers with greatest common divisor $1$.
Also in this case the conjecture turned out to be true. Many more special cases 
have been proved since. The reader is refered to the paper [EWZ], where
several of them have been established.

Studying the Generalized
Vanishing Conjecture, Zhao observed that
the main obstruction to understand the condition $\Lambda^m(P^m)=0$
is the fact that the differential operator $\Lambda$ and the
multiplication operator $P$ do not commute. Therefore he considered
the {\em left symbol map} $\mathfrak{L}$ from the Weyl algebra $A_n(\C)$
to the polynomial ring $\C[\zeta,z]$, which is the $\C$-linear map
sending each monomial $\partial^a z^b$ to $\zeta^a z^b$. So to compute
the image of a differential operator under this map, one first needs to write
the operator as a $\C$-linear combination of monomials of the form
$\partial^a z^b$ and then replace each $\partial_i$ by $\zeta_i$. 
In a similar way one can define the {\em right symbol map}
$\mathfrak{R}$ from the Weyl algebra to the polynomial ring
$\C[\zeta,z]$ by defining $\mathfrak{R}(z^a\partial^b)=z^a\zeta^b$.

Now let $\circ$ denote the multiplication in the Weyl algebra $A_n(\C)$,
i.e.\@ the composition as $\C$-linear maps acting on the polynomial ring
$\C[z]$. Then the condition $\Lambda^m(P^m)=0$ is equivalent to
$(\Lambda^m\circ P^m)(1)=0$. Furthermore, the differential operators which 
annihilate the element $1$ form the left ideal in $A_n(\C)$ generated
by the partial derivatives $\partial_i$. Applying the right symbol map
$\mathfrak{R}$, we therefore obtain that $\Lambda^m(P^m)=0$, if and only if
$\mathfrak{R}(\Lambda^m \circ P^m)$ belongs to the ideal generated
by the $\zeta_i$ in $\C[\zeta,z]$, or equivalently that
$\pi\circ\mathfrak{R}(\Lambda^m \circ P^m)=0$, where
$\pi$ denotes the $\C[z]$-homomorphism from $\C[\zeta,z] $ to $\C[z]$
sending each $\zeta_i$ to zero. Finally observe that

$$\mathfrak{L}(\Lambda^m\circ P^m)=\Lambda(\zeta)^mP(z)^m.$$

\noindent Since $\mathfrak{L}$ is an isomorphism of $\C$-vector spaces
this implies that

$$\Lambda^m\circ P^m=\mathfrak{L}^{-1}(\Lambda({\zeta})^m P(z)^m).$$

\noindent Combining this with the observation above, we obtain that

\medskip

$\Lambda^m(P^m)=0$ if and only if $L((\Lambda(\zeta)P(z))^m)=0$, 

\medskip

\noindent where $L=\pi\circ\mathfrak{R}\circ\mathfrak{L}^{-1}$.
In a similar way the condition $\Lambda(QP^m)=0$ is equivalent
to $L(Q(z)(\Lambda(\zeta)P(z))^m)=0$. So combining both results
above the Generalized Vanishing Conjecture can be reformulated
as follows: let $f=\Lambda(\zeta)P(z)$ and $g=Q(z)$. If $f^m$
belongs to $kerL$ for all positive $m$, then $gf^m$ belongs to
$kerL$ for almost all $m$. Having Mathieu's observations in mind
it then was a minor step to generalize this conjecture to
the stronger statement that the above implication should
hold for all $f$ and $g$ in $\C[\zeta,z]$. In other words 
these calculations led Zhao to conjecture that $kerL$ is
a Mathieu subspace of $\C[\zeta,z]$.

Then the next natural question to consider is: is there a nice
way to describe $kerL$? To answer this question, let's look
at the definition of $L$. It is an easy excercise to verify
that $L(\zeta^a z^b)=\partial^a(z^b)$. In other words the
image of a polynomial in $\zeta$ and $z$ under $L$ is obtained
as follows: write the polynomial as a $\C$-linear combination of 
monomials of the form $\zeta^a z^b$, then replace 
$\zeta$ by $\partial$ and apply $\partial$ to the monomial in $z$.
From this description, one deduces readily that for each polynomials
$g$ in $\C[\zeta,z]$, the element of $\partial_{z_i}(g)-\zeta_i g$ belongs
to $kerL$ and more generally that $ImD\subset kerL$, where
$D$ is the set of $n$ commuting operators $\partial_{z_i}-\zeta_i$ and
$ImD$ is as defined in section one. Then finally Zhao showed that in fact
$ImD$ is equal to $kerL$, which by the above arguments led to his
original formulation of the Image Conjecture. A proof of the
equality of $kerL$ and $ImD$ is given in the next section,
where we also show that, as indicated by the arguments
above, the Image Conjecture (for the ring $A=\C[\zeta_1,\ldots,\zeta_n]$)
implies the Generalized Vanishing Conjecture.

\section{The Image Conjecture and the Generalized Vanishing Conjecture}

In this section we show that the Image Conjecture implies the Generalized
Vanishing Conjecture (and hence the Jacobian Conjecture). More precisely,
we only need the following special case of the Image Conjecture:
$A=\C[\zeta]=\C[\zeta_1,\ldots,\zeta_n]$ and $a_i=\zeta_i$ for each $i$.
Let's denote this case by $IC(n)$ for simplicity.

\medskip

\noindent {\bf Theorem 3.1.} {\em $IC(n)$ implies $GVC(n)$.}

\medskip

\noindent Let $L$ be the $\C$-linear map from $\C[\zeta,z]$ to $\C[z]$
defined in the previous section. In other words

$$L(\zeta^{a}z^{b})={\partial_z}^a(z^b).$$

\noindent We will show below that we have the equality

$$ImD=kerL.$$

\noindent Before we prove this result let us first show how it
can be used to prove the theorem.
So let $\Lambda=\Lambda(\partial)$ be a differential
operator with constant coefficients and $P$ a polynomial in $\C[z]$
such that $\Lambda^m(P^m)=0$ for all positive $m$. Let $g\in\C[z]$.
We must show that also $\Lambda^m(gP^m)=0$ for almost all $m$.
Therefore put

$$f(\zeta,z)=\Lambda(\zeta)P(z).$$

\noindent It then follows from the hypothesis and the definition of
$L$ that $L(f^m)=0$ for all positive $m$. Using the
equality $kerL=ImD$, this implies that all positive powers
of $f$ belong to $ImD$, which by assumption is a Mathieu subspace
of $\C[\zeta,z]$. Consequently also $gf^m$ belongs to $ImD=kerL$
for almost all $m$. So $L(gf^m)=0$ for almost all $m$. Using the
definition of $f$ this gives that $L(\Lambda(\zeta)^m g(z)P(z)^m)=0$,
whence $\Lambda^m(gP^m)=0$ for almost all $m$. So $GVC(n)$ holds.

\medskip

To conclude this section we prove the equality $ImD=kerL$.
For another proof we refer the reader to [Z3]. Our proof
is based on the following result from the theory of $\mathfrak{D}$-modules

\medskip

\noindent {\bf Proposition 3.2.} {\em Let $M$ be an
$A_n(\C)=\C[t_1,\ldots,t_n,\partial_1,\ldots,\partial_n]$-module such
that each $\partial_i$ is locally nilpotent on $M$. Then each $f$ in $M$
can be written uniquely in the form $f=\sum t^a f_a$, where each $f_a$
belongs to $N=\bigcap_i ker(\partial_i,M)$. In particular $f\in\sum t_iM$
if and only if $f_0=0$.}

\medskip To obtain the equality $ImD=kerL$ we apply this proposition
to $M=\C[\zeta,z]$, which is an $A_n(\C)$-module by defining

$$t_if=(\zeta_i-\partial_{z_i})f \mbox{ and } \partial_i f=\partial_{\zeta_i}(f).$$

\noindent Then $N=\bigcap ker \partial_{\zeta_i}=\C[z]$. So $f\in\C[\zeta,z]$
can be written uniquely in the form

$$f=\sum (\zeta_i-\partial_{z_i})^a f_a(z)\,\,\,\,\,\,\,\,\,(*)$$

\noindent for some $f_a(z)$ in $\C[z]$, and in particular
$f\in\sum (\partial_{z_i}-\zeta_i)\C[\zeta,z]=ImD$ if $f_0(z)=0$.
Finally observe that (*) implies that $L(f)=f_0(z)$. So we get
that $f$ belongs to $ImD$ if $f$ belongs to $kerL$.
So $kerL\subset ImD$. Since we saw that $ImD\subset kerL$ in the previous section,
the desired equality follows. 

\medskip

\noindent The proof of the proposition above follows
easily by applying the next lemma repeatedly

\medskip

\noindent {\bf Lemma 3.3.} {\em Assume that $M$ is an $A_1(\C)=\C[t,\partial]$-module 
and $f$ an element of $M$ such that $\partial^m f=0$ for some $m$.
Then $f=f_0+tf_1+\ldots+t^{m-1}f_{m-1}$ for some $f_i\in N=ker\partial$
which are uniquely determined.}

\medskip

\noindent {\em Proof}. The uniqueness follows easily by applying $\p$
sufficiently many times. So let $\p^m f=0$. Then $\p^{m-1}(\p f)=0$,
so by induction on $m$ we get

$$\p f=g_0+tg_1+\ldots+t^{m-2}g_{m-2}$$

\noindent for some $g_i$ in $N$. Now let

$$G:=\int \partial f=tg_0+\frac{1}{2}t^2g_1+\ldots+\frac{1}{m-1}t^{m-1}g_{m-2}.$$

\noindent Then $\p G=\p f$, so $f-G\in N$. Say $f-G=f_0$
for some $f_0$ in $N$. Using the definition of $G$
the desired result follows.

\section{A surprising connection}

In the mean time Zhao had taken a completely different approach.
He just wondered if sequences of the form $\Lambda P,\Lambda^2(P^2),
\Lambda^3(P^3),\ldots$, where $P$ can be any polynomial, had been studied
before. At the same time he took a closer look at the Laplace operator by 
compactifying real $n$-space, which led him to investigate
eigen functions of this operator on the $n$-sphere. Studying
the literature he came in contact with the Gegenbauer polynomials,
a special class of orthogonal polynomials. In particular he found
the classical Rodriques' formula, which gives a useful way to
describe these polynomials. 
Before I continue this story, let me first recall some basic 
facts concerning orthogonal polynomials (see also [DX] and [S]).

\subsection*{Orthogonal polynomials}

\indent Let $B$ be a non-empty open subset of $\R^n$ and $w$ a so-called
{\em weight function} on $B$ i.e.,\@ it is strictly positive on $B$
and its integral over this set is finite and positive.
To such a function one can associate a {\em Hermitian inner product}
on the $n$-dimensional polynomial ring $\C[x]$ by defining

$$\langle f,g\rangle=\int_B f(x)\overline{g(x)}w(x)dx.$$

\noindent A set of polynomials $u_a$, where $a=(a_1,\ldots,a_n)$ runs through $\N^n$,
is called {\em orthogonal over $B$} with respect to the weight function
$w$, if they form an orthogonal basis of $\C[x]$ with respect to the
associated inner product described above and satisfy the additional condition 
that the degree of each polynomial $u_a$ is equal to $|a|$, the sum of all $a_i$.

A standard way to construct orthogonal polynomials in one variable
is to apply the Gram-Schmidt process to the basis $1,x,x^2,\ldots$
Making special choices for $B$ and $w$ gives the following
so-called {\em classical orthogonal polynomials}.\\
\noindent 1. The {\em Hermite polynomials}: $B=\R$ and $w(x)=e^{-x^{2}}$.\\
\noindent 2. The {\em Laguerre polynomials}: $B=(0,\infty)$ and
$w(x)=x^{\alpha}e^{-x}$, with $\alpha>-1$.\\
\noindent 3. The {\em Jacobi polynomials}: $B=(-1,1)$ and
$w(x)=(1-x)^{\alpha}(1+x)^{\beta}$, with $\alpha,\beta>-1$. In case
both parameters are zero, i.e.\@ $w=1$, the polynomials are called
{\em Legendre polynomials}.

    From univariate orthogonal polynomials one can construct orthogonal
polynomials in dimension $n$ as follows: for each $1\leq i\leq n$
choose on open subset $B_i$ of $\R$ and a weight function
$w_i$ on $B_i$. Let $u_{i,m}$ with $m\geq 0$ be an orthogonal set
of univariate polynomials with respect to $B_i$ and $w_i$.

Then $B=B_1\times\ldots\times B_n$ is an open subset
of $\R^n$ and $w$ defined by $w(x)=w_1(x_1)\ldots w_n(x_n)$ is
a weight function on $B$, where $x=(x_1,\ldots,x_n)$.
Furthermore one easily verifies that the polynomials

$$u_a(x)=u_{1,a_1}(x_1),\ldots,u_{n,a_n}(x_n)$$

\noindent where $a=(a_1,\ldots,a_n)$, form an orthogonal set
of polynomials over $B$ with respect to the weight function $w$.

The multivariate orthogonal polynomials obtained from Hermite polynomials
will again be called Hermite polynomials. Similarly we get multivariate
Laguerre and Jacobi polynomials. These polynomials we call the {\em classical
(multivariate) orthogonal polynomials}.

Now a surprising result is that all these classical orthogonal
polynomials can be obtained from the so-called Rodrigues' formula.
With the terminology introduced above it asserts the following

\medskip

\noindent {\bf Rodrigues' formula.} {\em Let $u_a$ be a system of classical
orthogonal polynomials. Then there exist a $g=(g_1,\ldots,g_n)$ in
$\C[x]^n$ and non-zero
real constants $c_a$ such that

$$u_a=c_aw^{-1}{\partial_x}^{|a|}(wg^a)$$}

\noindent For example in the one-dimensional case, one obtains the
Hermite, Laguerre and Jacobi polynomials by taking $g=1,x,1-x^{2}$
respectively and the constants $c_a$ are respectively equal to
$(-1)^a,\frac{1}{a!}$ and $\frac{(-1)^a}{2^aa!}$.

\indent Then Zhao made a remarkable discovery, namely if
one defines

$$\Lambda_i=w^{-1}\circ\partial_i\circ w=\partial_i+w^{-1}\partial_i(w)$$

\noindent and

$$\Lambda=(\Lambda_1,\ldots,\Lambda_n)$$

\noindent then 

$$u_a=c_a\Lambda^a(g^a)\,\,\,\,\,\,\,\,\,\,(1)$$

\noindent In other words all classical
orthogonal polynomials come from sequences of the form $\Lambda^a(g^a)$.
For example in the one-dimensional case, taking for $\Lambda$
the operators $\partial-2x,\partial+(\alpha x^{-1}-1)$ and
$\partial-\alpha (1-x)^{-1}+\beta (1+x)^{-1}$, one gets the Hermitian,
Laguerre and the Jacobi polynomials respectively, apart from
the constants described above.

\medskip

\noindent Zhao was struck by the fact that, just as in the formulation
of the Image Conjecture, again differential operators of order one appeared.
He therefore wondered if a similar kind of Image Conjecture would
hold for the commuting set of differential operators $\Lambda$ coming from 
the classical orthogonal polynomials as described above. He made the
following modified conjecture

\medskip

\noindent {\bf Image Conjecture for Classical Orthogonal Polynomials.}\\
\noindent {\em Let $\Lambda$ be as described above and let

$$Im^{'}\Lambda:=\C[z]\cap (\sum_i \Lambda_i(\C[z]))\,\,\,\,\,\,\,\,(2)$$

\noindent Then $Im^{'}\Lambda$ is a Mathieu subspace of $\C[z]$.}

\medskip

Now let's take a closer look at the intersection described
in (2). Using the notation for $g$ and $u_a$ introduced above, one can verify
by explicit calculation that in the univariate case 
$\Lambda^m(g^a)$ is a polynomial, if $m$ is at most $a$. From this and (1)
one obtains that each multivariate classical orthogonal polynomial $u_a$,
where $a$ is not the zero vector, belongs to $Im^{'}\Lambda$. These observations
lead to the following interesting result.

\medskip

\noindent {\bf Proposition 4.1.} {\em Notations as above. Assume that $1$ does
not belong to $Im^{'}\Lambda$. Then\\
\noindent i) $Im^{'}\Lambda$ is the $\C$-linear span of all $u_a$
where $a$ is non zero.\\
\noindent ii) $Im^{'}\Lambda=\{f\in\C[x]\,|\, \int_B fwdx=0\}$.}

\medskip

\noindent {\em Proof}. i) follows from (2 ), the hypothesis and the
observation above proposition 4.1.\\
\noindent ii) Let $f\in\C[x]$. Write it in the basis $\{u_a\}$,
say $f=\sum f_a u_a$ with $f_a$ in $\C$. Since $u_0\in\C^*$, it follows
from i) that $f$ belongs to $Im^{'}\Lambda$ if and only if $f_0=0$.
Since the $\{u_a\}$ form an orthogonal basis and $u_0\in\C^*$ we have
that

$$f_0\langle u_0,u_0\rangle=\langle f,u_0\rangle=\int_B f\cdot \overline{u_0}\cdot wdx=
\overline{u_0}\int_B fwdx.$$

\noindent So $f$ belongs to $Im^{'}\Lambda$ if and only if $\int_B fwdx=0$.

\medskip

\indent Inspired by Mathieu's conjecture, replacing the connected
compact Lie group $G$ by the open subset $B$ of $\R^n$, the Haar measure
by a positive measure $d\sigma$ and the $G$-finite functions by
polynomials, the above proposition led Zhao to the following analogue 
of Mathieu's conjecture

\medskip

\noindent {\bf Zhao's Integral Conjecture.} {\em Let $B$ be an open subset of $\R^n$ and
$\sigma$ a positive measure on $B$ such that for any polynomial
$g$ in $\C[x]$ the integral $\int_Bgd\sigma$ is finite. Then the set
of polynomials $f$ whose integral over $B$ is zero, is a Mathieu
subspace of $\C[x]$.}

\medskip

\noindent This conjecture is widely open in all dimensions. 
We will return to it in section 6 below. For now we show

\medskip

\noindent {\bf Corollary 4.2.} {\em The Integral Conjecture implies the Image
Conjecture for classical orthogonal polynomials.}

\medskip

\noindent {\em Proof.} Since a weight function $w$ is strictly positive on $B$,
the measure $d\sigma=wdx$ is positive on $B$. Consequently if $1$
does not belong to $Im'\Lambda$, the above
proposition together with the Integral Conjecture imply that this
set is a Mathieu subspace of $\C[x]$. In case $1$ does belong to
$Im'\Lambda$, the observation immediately before proposition 4.1 implies that 
this image is the whole polynomial ring and hence it is
a Mathieu subspace as well.

\medskip

\indent The Image Conjecture and the Image Conjecture for classical
orthogonal polynomials inspired Zhao in [Z3] to formulate the following
rather general statement

\medskip

\noindent {\bf General Image Conjecture.} {\em Let $k$ be a field, $A$
a commutative $k$-algebra and $A[z]=A[z_1,\ldots,z_n]$ the polynomial
ring over $A$. Let $D$ be a commuting set of differential operators
of the form

$$\sum_{i=1}^n c_i\partial_i +g(z)$$

\noindent where the $c_i$ belong to $A$ and $g(z)\in A[z]$. Then

$$ImD:=\sum_{\Lambda\in D} \Lambda A[z]$$

\noindent is a Mathieu subspace of $A[z]$.}

\medskip

\noindent {\bf Comment 1.} If the field $k$ has characteristic zero,
no counterexamples to this conjecture are known. On the other hand,
as we have seen in the example at the end of section one, the conjecture 
is false if the characteristic of $k$ is positive. At this moment it is not 
clear what extra condition, similar to the one given in the
statement of the Image Conjecture, can be added to avoid this kind of obvious
counterexamples.

\medskip

\noindent {\bf Comment 2.} It is shown in [Z3] that if the field $k$
has characteristic zero,
then in order to prove or disprove the General Image Conjecture 
for sets of operators for which the $c_i$ belong to $k$ (not just in $A$), it
suffices to study the cases where the set $D$ consists of $n$
operators of the form $\p_i-\p_i(q(z))$, where $q(z)$ is some polynomial.
Observe that if we take for $q(z)$ the linear form $a_1z_1+\ldots+a_nz_n$
we obtain the statement of the Image Conjecture as described in section one.

\medskip

\indent After having described various conjectures, it is time to
investigate the question: what evidence is there in favour
of these conjectures and in particular what evidence supports
the Image Conjecture?

\section{The Image Conjecture in positive characteristic}

Throughout this section (except in the crucial lemma below) $k$ 
will be a field of characteristic $p>0$. With the notations 
introduced in section one we get

\medskip

\noindent {\bf Theorem 5.1} {\em If $(a_1,\ldots,a_n)$ is a regular 
sequence in $A$, then $ImD$ is a Mathieu subspace of $A[z]$. In other
words the Image Conjecture is true for positive characteristic.}

\medskip

\noindent The proof of this theorem is based on the following
result, whose proof will be sketched at the end of this section
(we refer to [EWrZ] for more details).

\medskip

\noindent {\bf Crucial lemma.} {\em Let $k$ be any field. Let $b$ be 
a polynomial of degree $d$ and denote by $b_d$ it homogeneous 
component of degree $d$. If $b$ belongs to $ImD$, then all 
coefficients of $b_d$ belong to the ideal $I$ of $A$ generated 
by all $a_i$.}

\medskip

\noindent {\bf Corollary 5.2.} {\em Let $f$ be a sum of monomials $f_a z^a$.
If $f^p$ belongs to $ImD$ then each $f_a^p$ belongs to $I$.}

\medskip

\noindent {\em Proof}. Write $f$ as a sum of homogeneous components
$f_i$. Then $f^p$ is a sum of $f_i^p$. It then follows from the
crucial lemma that all coefficients of $f_d^p$ belong to $I$,
where $d$ is the degree of $f$. So $f_d^p$ is a sum of monomials
of the form $c a_i z^{ap}$, with $|a|=d$. Since each such a
monomial is equal to $(\p_i-a_i)(-cz^{ap})$, which belongs to $ImD$,
it follows that $f_d^p$ belongs to $ImD$. Substracting this
polynomial from $f^p$ we obtain that 

$$f_0^p+\ldots+f_{d-1}^p \in Im D.$$

\noindent Then the result follows by induction on $d$.

\medskip

\noindent {\bf Proof of the theorem 5.1.} \\
\noindent i) Since $\p^p=0$ on $A[z]$, we get
that $a_i^p z^a=(\p_i-a_i)^p(-z^a)\in ImD$. So 
a polynomial belongs to $ImD$ if all its coefficients
belong to the ideal $J$ generated by all $a_i^p$.\\
\noindent ii) Now let $f$ be such that $f^p$ belongs to 
$ImD$ and $g$ be any polynomial. By i) it suffices to show
that all coefficients of $gf^m$ belong to $J$ if $m\geq p^2$.
Therefore write $f$ as a sum of monomials $f_a z^a$. Since $f^p$ 
belongs to $ImD$, it follows from corollary 5.2 that each 
$f_a^p$ belongs to $I$ and hence each $f_a^{p^2}$ belongs to $J$. 
Since $f^{p^2}$ is a sum of the monomials $f_a^{p^2}z^{p^2a}$, 
it then readily follows that all coefficients of $gf^m$ belong
to $J$ if $m\geq p^2$. As observed above this concludes the proof.

\medskip

\noindent {\em Proof of crucial lemma (sketch).}

\medskip

\noindent We only sketch the case $n=2$. Let $b\in ImD$,
say 

$$b=(\p_1-a_1)p+(\p_2-a_2)q \,\,\,\,\,\,\,\,\,\,(*)$$

\noindent for some polynomials $p$ and $q$.
Let $d$ be the degree of $b$ and denote by $b_d$ the homogeneous
component of degree $d$. Now we assume for simplicity that the degrees
of both $p$ and $q$ are at most $d+2$. Then looking at the component
of degree $d+2$ in (*) we get

$$-a_1p_{d+2}-a_2q_{d+2}=0$$

\noindent where $p_i$ and $q_i$ denote the homogeneous components
of degree $i$ of $p$ and $q$ respectively. From the regularity hypothesis 
on the sequence $a_1,a_2$ it then follows that there exists a 
polynomial $g_{d+2}$, homogeneous of degree $d+2$, such that 

$$p_{d+2}=a_2g_{d+2} \mbox{ and } q_{d+2}=-a_1g_{d+2}.$$

\noindent Comparing the components of degree $d+1$ in the equation (*) we get

$$\p_1p_{d+2}+\p_2q_{d+2}-a_1p_{d+1}-a_2q_{d+1}=0.$$

\noindent Substituting the formulas for $p_{d+2}$ and $q_{d+2}$, found above, gives

$$a_1(p_{d+1}+\p_2g_{d+2})-a_2(q_{d+1}-\p_2g_{d+2})=0.$$

\noindent Again, from the regularity condition it then follows that there
exists some polynomial $g_{d+1}$, homogeneous of degree $d+1$, such that

$$p_{d+1}=-\p_2g_{d+2}+a_2g_{d+1} \mbox{ and } q_{d+1}=\p_1g_{d+2}+a_1g_{d+1}.$$

\noindent Comparing the components of degree $d$ in the equation (*) we get

$$b_d=\p_1p_{d+1}+\p_2q_{d+1}-a_1p_d-a_2q_d.$$

\noindent Finally, substituting the formulas for $p_{d+1}$ and $q_{d+1}$ in the 
last equality gives

$$b_d=-a_1(p_d-\p_2g_{d+1})-a_2(q_d-\p_1g_{d+1}).$$

\noindent which gives the desired result.

\section{Examples of Mathieu subspaces}

Before we give some more evidence in favour of the Image Conjecture
we want to discuss various examples of Mathieu spaces.
The first example concerns the ring of $n\times n$ matrices
over a field $k$, where either the characteristic of $k$ is zero or
greater than $n$.

\medskip

\noindent {\bf Example 6.1.} {\em Let $R$ be the ring of $n\times n$
matrices over $k$. Then the subspace consisting of all matrices
having trace zero is a Mathieu  subspace of $R$.}

\medskip

\noindent Indeed, it is well-known that if the traces of the first
$n$ powers of a matrix $A$ are zero, then the matrix is nilpotent
and hence its $n$-th power is the zero matrix. Consequently
for any matrix $B$ also $BA^m=0$ if $m$ is at least $n$.
In particular the trace of this matrix is zero for all $m\geq n$.

\medskip

\noindent {\bf Remark 6.2.} {\em It is proved in [Z5] that, 
under the hypthesis on the characteristic of $k$ described above, 
the subspace of Example 6.1 is the only co-dimension one
left or right Mathieu subspace of $R$. In case the
characteristic is positive and at most $n$, it turns out that $R$
has no left or right Mathieu subspaces of co-dimension one.}

\medskip

\indent Now let me give two less trivial examples. Both concern
subspaces of the univariate polynomial ring $\C[t]$. Again
the reader is invited to find elementary proofs for the 
following two results.

\medskip

\noindent {\bf Example 6.3.} {\em The set of all polynomials in $\C[t]$
such that $\int_0^1 f dt=0$, is a Mathieu subspace
of $\C[t]$. In fact, if a polynomial $f$ is such that $\int_0^1 f^m dt=0$
for almost all $m$, then $f=0$.}

\medskip

\noindent {\bf Example 6.4.} {\em The set of all polynomials $f$ in $\C[t]$
such that $\int_0^{\infty} f e^{-t} dt=0$, is a Mathieu space 
of $\C[t]$. In fact, if a polynomial $f$ is such that 
$\int_0^{\infty} f^m e^{-t}dt=0$ for almost all $m$, then $f=0$.}

\medskip

\noindent A beautiful proof of the result described in Example 6.3 was given
by Mitya Boyarchenko in a personal communication to Zhao. This proof
has been included in the recent preprint [FPYZ]. Using his
techniques we were able to prove the statement described in Example 6.4 (see [EWrZ]
and section 8 below).
The importance of this result is due to the fact that it implies the one dimensional
Image Conjecture $IC(1)$ (for a proof we refer to the next section).
This implication reveals a remarkable fact, namely that the classical
Laguerre polynomials play a special role in the study of the Image Conjecture and
hence also of the Jacobian Conjecture.

To conclude this section let us observe
that the apparently stronger statements in the second half of both
Examples 6.3 and 6.4, are in fact equivalent to the ones made in the first
halves of these examples.
To see this, we need the following result of ([EZ2]).

\medskip

Let $B$ be an open set of $\R$ and $\sigma$ a positive
measure on $\R$ such that $\int_B f d\sigma$ is finite for all $f\in\C[t]$ and
such that the $\C$-bilinear form defined by $\langle f, g\rangle=\int_B fg d\sigma$
is non-singular, i.e.\@ for each non-zero $f$ there exists a $g$ such that
$\langle f, g\rangle$ is non-zero.

\medskip

\noindent {\bf Proposition 6.5.} {\em If the set of all
polynomials $f$ with $\int_B f d\sigma$ is a Mathieu subspace
of $\C[t]$, then the only polynomial $f$ such that
$\int_B f^m d\sigma=0$ for almost all $m$, is the zero polynomial, i.e.\@ $f=0$.}

\medskip

\indent The announced equivalence of the
statements made in both Examples 6.3 and 6.4 then follows by taking
$B=(0,1),d\sigma=dt$ and $B=(0,\infty),d\sigma=e^{-t}dt$ respectively: in both cases
the hypothesis of the proposition is satisfied, since the corresponding bilinear
forms are in fact Hermitian inner products on $\C[t]$.

\section{The one dimensional Image Conjecture}

In this section we show how the result of Example 6.4 implies $IC(1)$.

\medskip

\noindent {\bf Theorem 7.1.} {$IC(1)$ is true.}

\medskip

\noindent {\em Proof.} Let $L$ be the $\C$-linear map from $A=\C[\zeta,z]$
to $\C[z]$ defined in section 2 and 3 by the formula $L(\z^az^b)=\p^a(z^b)$.
So this expression is zero if $a$ is larger than $b$. Furthermore for any non-zero polynomial $g$ we define its degree, denoted $\Deg(g)$, as the maximum of the degrees of all non-zero monomials appearing in $g$, where $\Deg(c\z^a z^b)=b-a$ and $c$ is a non-zero constant in $\C$. It follows that

$$\mbox{ if } \Deg(g)\leq -1, \mbox{ then } L(g)=0\,\,\,\,\,\,\,\,(3)$$

\noindent In particular, if for some element $f$ of $A$ its degree
is $\leq -1$, then certainly the degrees of all powers $f^m$ are $\leq -1$,
which by (3) implies that $L(f^m)=0$ for all $m\geq 1$. Now we will
show that the converse is true as well. 

\medskip

\noindent {\bf Proposition 7.2.} Let $f\in A$. Then $L(f^m)=0$ for 
almost all $m$, if and only if $\Deg(f)\leq -1$.

\medskip

\noindent Before we prove this proposition let us show
how it implies Theorem 7.1. So assume that $L(f^m)=0$ for almost all $m$. Then by the proposition $\Deg(f)\leq -1$, so the degree of $f^m$ is at most $-m$. Now let $g$ in $A$ be non-zero and let $d$
be its degree. It then follows that the degree of $gf^m$ is at most $-1$,
if $m$ is at least $d+1$. The result then follows from (3).

\medskip

\indent To prove proposition 7.2, let $f$ be such that $L(f^m)=0$ for almost
all $m$. {\em Assume} that $r:=\Deg(f)\geq 0$. We will arrive at a contradiction.
Namely, let $g=\z^r f$. Then $\Deg(g)=0$. Furthermore 

$$L(g^m)=L(\z^{mr}f^m)=\p^{mr}L(f^m)=0 \,\,\,\,\,\,\,\,\,(4)$$ 

\noindent for almost all $m$. Writing $g$ in its homogeneous decomposition, 
using that the degree of $g$ is zero,
we get that $g=g_0+g_{-1}+\dots$ and hence that $g^m=g_0^m+g_*$, where $g_0$ is non-zero
and $g_*$ has degree at most $-1$. Applying $L$ to the last equality it follows
from (3) and (4) that $L(g_0^m)=L(g^m)=0$ for almost all $m$. Summarizing,
if $r\geq 0$ there exists a non-zero element $g_0$, which is homogeneous
of degree zero, such that $L(g_0^m)=0$ for almost all $m$. 

Write again $g$
instead of $g_0$. Then $g$ is a sum of monomials of the form $c_a\z^az^a$. In
other words $g=P(u)$, a non-zero polynomial in $u:=\z z$ over $\C$. Now observe
that $L(u^n)=L(\z^n z^n)=n!$. Since, as one easily verifies, 
$\int_0^{\infty} u^n e^{-u} du=n!$ for all $n\geq 0$, it follows
that $L(f(u))=\int_0^{\infty} f(u) e^{-u} du$ for each polynomial $f(u)$.
Since $L(g^m)=L(P(u)^m)=0$ for almost all $m$, it then follows from
Example 6.4 that $P(u)=0$, i.e.\@ $g=0$, a contradiction since $g=g_0$ is non-zero.

\medskip

The proof given above shows that it is interesting and necessary to investigate $IC(n)$ for polynomials of the form $f(u_1,\ldots,u_n)$, where $u_i=\zeta_i z_i$.
One easily verifies that $L(u_1^{a_1}\cdots u_n^{a_n})=a_1!\cdots a_n!$.
This leads to the following question: if $L(f^m)=0$ for almost all $m$, does this imply that $f=0$?
We don't know the answer to this question when $n$ is at least two. However
various computations suggest that the following is true: if $f$ is a sum
of $N$ monomials and $L(f^m)=0$ for the first $N$ exponents $m$, then $f=0$.
Since, similar as in the proof of $IC(1)$ given above, one has the equality

$$L(f(u_1,\ldots,u_n))=\int_{(0,\infty)^n} f(u_1,\ldots,u_n) e^{-(u_1+\ldots+u_n)}du_1\ldots du_n$$
 
\noindent we make the following conjecture

\medskip

\noindent {\bf Conjecture 7.3.} {Let $\C[u]=\C[u_1,\dots,u_n]$ and $f\in\C[u]$.\\
\noindent  Let $B=(0,\infty)^n$. If $\int_B f(u)^m e^{-(u_1+\dots +u_n)}du=0$ 
for almost all $m$, then $f=0$.}

\medskip

\noindent Observe again, that this conjecture is a stronger version of a special 
case of Zhao's Integral conjecture. We refer the reader to the paper [EWrZ],
where several special cases of this conjecture are proved.

On the other hand, proving conjecture 7.3 is not enough
to prove the Image Conjecture $IC(n)$ in higher dimension. Namely the following result of [Z3] shows, that there exist many polynomials $f$ in $\C[\zeta,z]$,
which do not belong to $\C[u_1,\ldots,u_n]$, but do have the
property that $L(f^m)=0$ for all positive $m$.

\medskip

\noindent {\bf Proposition 7.4.} {\em Let $H=(H_1,\ldots,H_n)$ in $\C[z]^n$ be 
such that each $H_i$ has no terms of degree at most one. Let

$$f(\zeta,z)=\z_1 H_1+\ldots+\z_n H_n.$$

\noindent Then $L(f^m)=0$ for all positive $m$ if and only if
the Jacobian matrix of $H$ is nilpotent.}

\section{Final remarks on the Image Conjecture.}

To conclude this paper we will discuss a slightly stronger version
of $IC(n)$, which we denote by $IC_*(n)$. At this moment it is not known if it is really stronger. First we introduce some notations.\\
Let $\mathfrak{E}$ be a collection of fields. We make the following conjecture

\medskip

\noindent {\bf IC$_*(\mathfrak{E},n)$.} {\em For every pair of positive integers
$d$ and $e$ there exists a positive integer $D(d,e)$, such that the following
holds: for any field $k$ in $\mathfrak{E}$, any $b$ in $\N^n$ with $|b|=e$ and any polynomial $f$ in $\z,z$ over $k$ of degree $d$, such that $L(f^m)=0$, for 
all positive $m$, we have that $L(z^bf^m)=0$ for all $m\geq D(d,|b|)$.}

\medskip

\indent If $\mathfrak{E}$ consists of only one field $k$, we simply write
$IC_*(k,n)$ instead of $IC_*(\{k\},n)$. Finally, denote the set of
number fields by $\mathfrak{N}$.

\medskip

\noindent {\bf Reduction theorem.} {\em If $IC_*(\mathfrak{N},n)$ holds, then $IC_*(\C,n)$ holds.}

\medskip

\noindent {\em Proof.} i) First we show that the hypothesis implies that 
$IC_*(\overline{\Q},n)$ holds, where $\overline{\Q}$ is the algebraic closure
of $\Q$. Namely, let $f\in\overline{\Q}[\z,z]$, of degree $d$, be such that 
$L(f^m)=0$, for all positive $m$. Since $f$ has only a finite number 
of coefficients, which are all algebraic over $\Q$, it follows that 
all these coefficients belong to some number field. Then the hypothesis 
implies that, for each monomial $z^b$, there exists a positive integer
$D(d,|b|)$, such that  $L(z^bf^m)=0$ for all $m\geq D(d,|b|)$. 
So $IC_*(\overline{\Q},n)$ holds.\\
\noindent ii) Now we will show that $IC_*(\C,n)$ holds.
Let $d$ be some positive integer and denote 
by $f_U$ the universal polynomial of degree $d$ in $\z$ and $z$, i.e.\@ the coefficient
of each monomial $\zeta^pz^q$, where $|p|+|q|\leq d$, is a new variable $C_{p,q}$.
Denote by $\Z[C]$ the polynomial ring in these variables over $\Z$ 
and let $N$ denote the number of these variables.
For each positive $m$ we have that $L(f_U^m)$ is a polynomial in $z$, with
coefficients in $\Z[C]$. Let $I$ be the ideal in $\overline{\Q}[C]$ generated by all
these coefficients, i.e. for all positive $m$. By i) we know that $IC_*(\overline{\Q},n)$ holds. So for each monomial $z^b$ there exists
a positive integer $D(d,|b|)$ having the property decribed in i). 
For each $m\geq D(d,|b|)$ we get polynomials $L(z^bf_U^m)$ in $z$, with 
coefficients in $\Z[C]$.
The ideal in $\overline{\Q}[C]$ generated by all these coefficients we denote 
by $J_b$. 
By i) it follows that, if $c\in\overline{\Q}^N$ is a zero of $I$ and $z^b$
is some monomial, then $c$ is also a zero of $J_b$. Let $g_1,\ldots,g_s$
 be generators of the
ideal $J_b$. It then follows from the Nullstellensatz, that there exists some
natural number $r$, such that each polynomial $g_i^r$ belongs to $I$.

Now we can finish the proof. Namely, let $f$ be a polynomial of degree $d$ in $\z, z$
over $\C$, such that $L(f^m)=0$, for all positive $m$. Fix some monomial $z^b$. Let $c_{p,q}$ be the coefficient of $\z^p z^q$ in $f$. Then the vector $c$ in $\C^N$, whose components are formed by the $c_{p,q}$, is a zero of $I$ and hence of the $g_i$ (since each $g_i^r$ belongs to $I$). So $c$ is a zero of $J_b$, which means
that $L(z^bf^m)=0$, for all $m\geq D(d,|b|)$.

\subsection*{A proof of $IC_*(\C,1)$}

By the reduction theorem it suffices to prove $IC_*(k,1)$ for number fields $k$.
So let $k$ be such a field. Let $d$ be a natural number and $z^b$ some monomial
in the single variable $z$. Put $D(d,b)=b+1$. Suppose now that $f\in k[\z,z]$
is such that $L(f^m)=0$ for all positive $m$. Looking at the proof of
proposition 7.2, we see that it suffices to prove the following result

\medskip

\noindent {\bf Lemma 8.1.} {\em Let $k$ be a number field and $g$ in $k[u]$ such
that $L(g^m)=0$ for almost all $m$. Then $g=0$.}

\medskip

\noindent Here $L$ is the $k$-linear map from $k[u]$ to $k$
defined by $L(u^i)=i!$, for all $i$.

\medskip 

\noindent {\em Proof of lemma 8.1.} Assume that $g$ is non-zero. Then we may 
assume that $g=u^s+c_{s+1}u^{s+1}+\ldots+c_du^d$ for some $c_i$ in $k$.
Since $k$ is a number field, there exists, for almost all prime numbers $p$,
a non-archimedean valuation $v$ on $k$, such that $v(p)=1$ and $v(c_i)\geq 0$
for all $i$. Now choose such a $p$ large enough, with $L(g^p)=0$.
We claim that the equation $L(g^p)=0$ leads to a contradiction, which shows 
that our assumption, that $g$ is non-zero, is false.
To obtain this contradiction, first observe that

$$g^p=u^{sp}+\sum_{i=s+1}^d c_i^pu^{ip}+p\sum_{i=sp+1}^{dp-1}h_i(c)u^i$$

\noindent where $h_i(c)$ belongs to the subring of $k$, generated by the $c_i$
over $\Z$. In particular $v(h_i(c))\geq 0$ for all $i$. Now applying $L$
to the equality above, using that $L(g^p)=0$ and $L(u^i)=i!$ for all $i$, gives

$$0=(sp)!+\sum_{i=s+1}^d c_i^p(ip)!+p\sum_{i=sp+1}^{dp-1}h_i(c)i!$$

\noindent Observe that, if $i\geq s+1$, then $(ip)!=ip(sp)!n_i$, for some
natural number $n_i$ and that, if $i\geq sp+1$, then $i!=q_i(sp)!$, for
some natural number $q_i$.
Then, dividing the last equation by $(sp)!$, gives

$$1+\sum_{i=s+1}^d c_i^p(ip)n_i+p\sum_{i=sp+1}^{dp-1}h_i(c)q_i=0.$$

\noindent Finally observe that, each term in each of the two sums, has
a positive valuation at $v$, since $v(p)=1>0$ and both $v(c_i),v(h_i(c))\geq 0$,
for all $i$. Since $v(1)=0$, this gives a contradiction.

\medskip

\noindent {\bf Acknowledgments.} This paper is an extended version of lectures
given at the University of Washington, St. Louis. The author likes to thank
this university for its great hospitality during his stay and in particular
David Wright for many stimulating conversations and being a wonderful host.
Ofcourse, a special word of thanks to Wenhua Zhao, for coming up with
these fascinating conjectures and for having various conversations on many of his 
wonderful ideas and results, both at the Illinois State University, at Normal
and at the Radboud University of Nijmegen. Finally the author is grateful to
Stefan Maubach, for inviting Wenhua to Nijmegen in April-May 2010 and to Michiel
de Bondt for his careful reading of the manuscript.

\section*{References}

\noindent [BCW] H. Bass, E. Connell, D. Wright, The Jacobian Conjecture: Reduction
of Degree and Formal Expansion of the Inverse, {\em Bull. Amer. Math. Soc.(NS)} 7 (1982), 287-330.\\
\noindent [B] J.-E. Bj\"ork, Rings of differential operators, {\em North-Holland
Math. Library Ser.}, Vol. 21, North-Holland, 1979.\\
\noindent [BE] M. de Bondt and A. van den Essen, A Reduction of the Jacobian
Conjecture to the Symmetric Case, {\em Proc. Amer. Math. Soc.} 133 (2005), no.8,
2201-2205.\\
\noindent [DvK] J.J. Duistermaat and W. van der Kallen, Constant terms in powers
of a Laurent polynomial, {\em Indag. Math.(NS)} 9 (1998), no.2, 221-231.\\
\noindent [DX] C. Dunkl and Y. Xu, Orthogonal polynomials in several variables,
Encyclopedia of Mathematics and its Applications, 81. {\em Cambridge Univ.
Press}, Cambridge, 2001.\\
\noindent [E] A. van den Essen, Polynomial Automorphisms and the Jacobian
Conjecture, {\em Prog. Math.}, Vol.190, Birkh\"auser Verlag, Basel, 2000.\\
\noindent [EWZ] A. van den Essen, R. Willems and W. Zhao, Some results on the \allowbreak
vanishing conjecture of differential operators with constant coefficients, \allowbreak arXiv: 0903.1478v1.\\
\noindent [EWrZ] A. van den Essen, D. Wright and W. Zhao, On the Image
Conjecture, In preparation.\\
\noindent [EZ] A. van den Essen and W. Zhao, Two results on homogeneous Hessian nilpotent
polynomials, {\em J. Pure Appl. Algebra} 212 (2008), no.10, 2190-2193.\\
\noindent [EZ2] A. van den Essen and W. Zhao, The Vanishing conjecture on
univariate polynomial rings. In preparation.\\
\noindent [FPYZ] J.P. Francoise, F. Pakovich, Y. Yomdin and W. Zhao, Moment
Vanishing Problem and Positivity; Some Examples. Preprint 2010.\\
\noindent [K] O.H. Keller, Ganze Cremona-Transformationen, 
{\em Monatsh. Math. Phys.}, 47 (1939), 299-306.\\
\noindent [M] O. Mathieu, Some conjectures about invariant theory and their
applications, Alg\`ebre non commutative, groupes quantiques et invariants (Reims, 1995),
S\'emin. Congr. 2, Soc. Math. France, 263-279 (1997).\\
\noindent [S] G. Szeg\"o, Orthogonal Polynomials, 4th edition, {\em Amer. Math. Soc., Colloq. Publ.}, 1975.\\
\noindent [Y] A. Yagzhev, On Keller's problem, {\em Siberian Math. J.}, 21 (1980),
747-754.\\
\noindent [Z1] W. Zhao, Hessian nilpotent polynomials and the Jacobian
Conjecture, {\em Trans. Amer. Math. Soc.}, 359 (2007), no.1, 294-274.\\
\noindent [Z2],--,A Vanishing Conjecture on differential operators with constant
oefficients, {\em Acta Math. Vietnam.}, Vol.32 (2007), no.3, 259-286.
See also arXiv:0704.1691v2.\\
\noindent [Z3],--,Images of commuting differential operators of order one with
constant leading coefficients, {\em J. Algebra} 324 (2010), no.2, 231-247.\\
\noindent [Z4],--, Generalizations of the Image Conjecture and the Mathieu
Conjecture, {\em J. Pure Appl. Algebra}, 214 (2010), 1200-1216.\\
\noindent [Z5],--, Mathieu subspaces of associative algebras, arXiv:1005.4260v1.\\

\medskip

\noindent Department of Mathematics, Radboud University Nijmegen\\
\noindent Postbus 9010, 6500 GL Nijmegen, The Netherlands\\
\noindent Email: essen@math.ru.nl

\end{document}